\documentclass[letterpaper,12pt]{article}

\usepackage{color}
\usepackage{amssymb}
\usepackage{amsmath}
\usepackage{amscd}
\usepackage{mathtools}
\usepackage[all]{xy}
\usepackage{tikz-cd}

\usepackage{hyperref}
\hypersetup{colorlinks=true,linkcolor=blue,anchorcolor=blue,citecolor=blue}
\usepackage{cleveref}

\usepackage{color}

\newtheorem{theo}{Theorem}[section]
\newtheorem{prop}[theo]{Proposition}
\newtheorem{lem}[theo]{Lemma}
\newtheorem{cor}[theo]{Corollary}
\newtheorem{defi}[theo]{Definition}
\newtheorem{rem}[theo]{Remark}

\def \id {{\rm id}}
\def \Br {{\rm{Br}}}

\def \Ga {{\Gamma}}

\def \Pic {{\rm {Pic}}}
\def \Div {{\rm {Div}}}
\def \div {{\rm{div}}}
\def \Gal {{\rm{Gal}}}
\def \Ker {{\rm{Ker}}}
\def \Im {{\rm {Im}}}

\def \A{{\mathbb A}}
\def \P{{\mathbb P}}

\def \dim {{\rm{dim}}}
\def \Hom {{\rm {Hom}}}

\def \Z {{\mathbb Z}}

\def \Ext {{\rm Ext}}

\def \s {{\rm s}}

\def\lra{\longrightarrow}

\def\H{{\rm H}}

\def\Ga{\Gamma}

\def\et{{\rm{\acute et}}}
\def\nr{{\rm nr}}

\newcommand{\bthe}{\begin{theo}}
\newcommand{\ble}{\begin{lem}}
\newcommand{\bpr}{\begin{prop}}
\newcommand{\bco}{\begin{cor}}
\newcommand{\bde}{\begin{defi}}
\newcommand{\brem}{\begin{rem}}
\newcommand{\ethe}{\end{theo}}
\newcommand{\ele}{\end{lem}}
\newcommand{\epr}{\end{prop}}
\newcommand{\eco}{\end{cor}}
\newcommand{\ede}{\end{defi}}
\newcommand{\erem}{\end{rem}}

\title{Low degree unramified cohomology of generic diagonal hypersurfaces}
\author{J.-L. Colliot-Th\'el\`ene and A.N.~Skorobogatov}
\date{\today}
\begin{document}
\maketitle

\begin{abstract}
\noindent We prove that the $i$-th unramified cohomology group of the generic diagonal hypersurface
in the projective space of dimension $n\geq i+1$ is trivial for $i\leq 3$.
\end{abstract}

\section{Introduction}

Let $k$ be a field with separable closure $k_\s$ and absolute Galois group $\Ga= \Gal(k_\s/k)$.
Let $\mu$ be a finite \'etale commutative group $k$-scheme of order not divisible by ${\rm char}(k)$.
The datum of such a group $k$-scheme
$\mu$ is equivalent to the datum of the finite $\Ga$-module $\mu(k_\s)$. 
For an integer $m\geq 2$ we denote by $\mu_m$ the group $k$-scheme of $m$-th roots of unity.
If $N$ is a positive integer not divisible by ${\rm char}(k)$ such that $N\mu=0$, then
$\mu(-1)$ denotes the commutative group
$k$-scheme ${\bf Hom}_{k-{\rm gps}}(\mu_{N},\mu)$. In terms of Galois modules,
$\mu(-1)$ is $\Hom_{\Z}(\mu_{N}(k_\s),\mu(k_\s))$ with the natural Galois action.

Let $X$  be a smooth integral variety over $k$.
We denote by $X^{(n)}$ the set of points of $X$ of codimension $n$.
In this paper, 
the unramified cohomology group $\H^i_\nr(X, \mu)$, where $i$ is a positive integer, is defined as
the intersection of kernels of the residue maps
$$\partial_x\colon\H^i(k(X), \mu)\to \H^{i-1}(k(x), \mu(-1)),$$
for all $x\in X^{(1)}$. 
For equivalent definitions, see \cite[Thm.~4.1.1]{CT95}. 
Restriction to the generic point of $X$ gives rise to a natural map
$$\H^i_\et(X,\mu) \to \H^i_\nr(X, \mu).$$
Purity for \'etale cohomology implies that it
is an isomorphism for $i=1$ and surjective for $i=2$, see \cite[\S 3.4]{CT95}.
In the case $i=2$ with $\mu=\mu_m$,
where $m$ is not divisible by ${\rm char}(k)$, this gives a canonical isomorphism
$$\Br(X)[m]\tilde\lra\H^2_\nr(X, \mu_m),$$
see \cite[Prop.~4.2.1 (a), Prop.~4.2.3 (a)]{CT95}.
If $X/k$ is smooth, proper, and integral,
then  $\H^i_\nr(X, \mu)$ does not depend on 
the choice of $X$ in its birational equivalence class, 
see \cite[Prop.~4.1.5]{CT95} and \cite[Remark (5.2), Cor.~(12.10)]{R96}.

Let  $n \geq 2$ and let $K=k(a_1,\ldots,a_n)$ 
be the field of rational functions in the variables $a_1,\ldots,a_n$.
Let $X_K\subset\P^n_K$ be
the hypersurface with equation
$$x_0^d+a_1x_1^d+\ldots+a_n x_n^d=0,$$
where $d$ is not divisible by ${\rm char}(k)$.
 In this paper,   for $i=1,2,3$ and $n\geq i+1$,
we prove that the natural map  
$$\H^i(K,\mu)\to\H^i_\nr(X_K,\mu)$$
is an isomorphism, see Theorem \ref{main}. 
In the case when $i=2$ and $\mu=\mu_m$ with $m\geq 2$, this gives that
the natural map of Brauer groups $\Br(K)\to\Br(X_K)$ is an isomorphism
of subgroups of elements of order not divisible by ${\rm char}(k)$, see Corollary \ref{Br}.
In the case when $k$ has characteristic zero, this result
was obtained as \cite[Thm.~1.5]{GS} by a completely different method,
using results on the topology of the Fermat surface as a complex manifold.

In this paper we use the  formalism proposed by M. Rost in \cite{R96} 
which applies {\em inter alia} to Galois cohomology  \cite[Remarks (1.11), (2.5)]{R96}.
We do not use the Gersten conjecture for \'etale cohomology \cite{BO74}.

Let us describe the structure of this note.
In Section \ref{unramifcoho} we recall some basic facts about unramified cohomology
including a functoriality property of the Bloch--Ogus complex 
with respect to faithfully flat morphisms with integral fibres.
In Section \ref{unram2hyp} we show that for smooth complete intersections $X\subset\P^n_k$
there are canonical isomorphisms $\H^i(k,\mu)\stackrel{\sim}\lra\H^i_\nr(X,\mu)$ for $i=1, 2$
when $\dim(X)\geq i+1$. Generic diagonal hypersurfaces are studied in Section \ref{gen diag}.
The easy proof of the main theorem in the case $i=1$ is given
in Section \ref{cohodegree1}.
This is used in the proof in the case $i=2$ in Section \ref{cohodegree2and3}, after
some preparations in Section \ref{basicdiagram}.
Finally, in Section \ref{pairs} we use a similar idea 
to give a short proof of the triviality of the Brauer group
of certain surfaces in $\P^3_{k(t)}$ defined by a pair of polynomials with coefficients in $k$,
see Theorem \ref{f=ag}, which was proved in \cite{GS} when ${\rm char}(k)=0$.

Our proof in this note develops a geometric idea suggested 
by Mathieu Florence during the second author's talk 
at the seminar ``Vari\'et\'es rationnelles" in November 2022. 
The authors are very grateful to Mathieu Florence for his suggestion.

\section{Functoriality of the Bloch--Ogus complex} \label{unramifcoho}

 For any smooth integral variety $X$ over $k$ and any $i\geq 2$ there is a complex
$$0\lra\H^{i}(k(X),\mu) \stackrel{(\partial_x)}\lra \bigoplus_{x \in X^{(1)}} \H^{i-1}(k(x),\mu(-1)) 
\stackrel{(\partial_y)}\lra  \bigoplus_{y\in X^{(2)}} \H^{i-2}(k(y), \mu(-2)),$$
which we call the {\em Bloch--Ogus complex}.
The maps in this complex are defined in \cite[(2.1.0)]{R96}.
(The map $\partial_{x}$ is the residue defined for discrete valuation rings
by Serre \cite{S03}, see also \cite[Def.~1.4.3]{CTS21}.)
The proof that the resulting sequence is a complex is given in \cite[Section 2]{R96}.
If $y \in X^{(2)}$ is a regular point
of the closure of $x \in X^{(1)}$, then the map
$\partial_y\colon\H^{i-1}(k(x),\mu(-1)) \to   \H^{i-2}(k(y), \mu(-2))$
is the residue map for the local ring of $y$ in the closure of $x$, which is a discrete valuation ring.

The unramified cohomology group $\H^i_\nr(X, \mu)$ is the
homology group of this complex at the term $\H^{i}(k(X),\mu)$,
i.e., the intersection of $\Ker(\partial_x)$ for all $x \in X^{(1)}$.
\medskip

Let $p \colon X \to Y$ be a faithfully  flat morphism of smooth integral $k$-varieties
with integral fibres. By \cite[Section (3.5); Prop.~(4.6)(2)]{R96},
there is a chain map of complexes 
$$\xymatrix{0\ar[r]&\H^{i}(k(X),\mu)\ar[r]& \bigoplus_{v \in X^{(1)}} \H^{i-1}(k(x),\mu(-1))   \ar[r]&  \bigoplus_{x \in X^{(2)}} \H^{i-2}(k(x), \mu(-2))      
\\0\ar[r]&\H^{i}(k(Y),\mu)    \ar[r]\ar[u]&\bigoplus_{y \in Y^{(1)}}
 \H^{i-1}(k(y),\mu(-1))       \ar[r]\ar[u]&\bigoplus_{y\in Y^{(2)}} 
\H^{i-2}(k(y), \mu(-2))      \ar[u]}$$
The middle vertical map is the natural one if $p(x)=y$, otherwise it is zero, and similarly
for the right-hand vertical map.

The morphism $X \to Y$ is called an {\em affine bundle} if Zariski locally on $Y$, 
it is isomorphic to $Y \times_{k} \A^n \to Y$ with affine transition morphisms.
In this case the vertical maps in the  above diagram induce isomorphisms on the  left-hand
and middle homology groups, see \cite[Prop.~(8.6)]{R96}. In particular, we have an isomorphism
\begin{equation}
\H^i_\nr(X, \mu )\cong\H^i_\nr(Y, \mu).\label{affine_bundle}
\end{equation}
Combined with  \cite[Cor.~(12.10)]{R96}, this implies that
$\H^i_\nr(X, \mu )$ is a stable birational invariant of smooth and proper 
integral $k$-varieties.

\section{Low degree unramified cohomology of complete intersections}\label{unram2hyp}

For a variety $X$ over a field $k$ we write $X^\s=X\times_k k_\s$.

\bpr\label{raetsel2}
 Let $X$ be a smooth, projective, geometrically integral variety over a field $k$
such that the natural map $\Pic(X)\to\Pic(X^\s)$ is an isomorphism of
 finitely generated free abelian groups.
 Then for any $k$-group of multiplicative type $M$ 
the natural map $$\H^2(k,M) \to \H^2(k(X),M)$$ is injective.
\epr
{\em Proof.} We have a commutative diagram with exact rows and natural vertical maps
\begin{equation}\begin{split}\xymatrix{
0\ar[r]& k_\s^\times\ar[r]& k_\s(X)^\times\ar[r]& \Div(X^\s)\ar[r]&
\Pic(X^\s)\ar[r]& 0\\
0\ar[r]& k^\times\ar[r]\ar[u]& k(X)^\times\ar[r]\ar[u]& \Div(X)\ar[r]\ar[u]&
\Pic(X)\ar[r]\ar[u]^\cong& 0}\label{el}
\end{split}\end{equation}
The abelian group $\Pic(X)$ is free, so the homomorphism $\Div(X)\to\Pic(X)$ has a section. 
Then our assumption implies that the map of $\Ga$-modules 
$\Div(X^\s)\to \Pic(X^\s)$ has a section.
By definition, the elementary obstruction 
$e(X)\in\Ext^2_k(\Pic(X_{\bar k}),\bar k^\times)$ is the class of the 2-extension of $\Ga$-modules
given by the upper row of (\ref{el}). Thus we have $e(X)=0$.
The result now follows from \cite[Prop.~2.2.5]{CS}. \hfill $\Box$

\ble \label{Lefschetz}
Let $X \subset \P^n_{k}$ be a complete intersection.

{\rm (a)} If $\dim(X)\geq 2$, then the natural map $\H^1(k,\mu)\to\H^1_\et(X,\mu)$
is an isomorphism. 

{\rm (b)} If $\dim(X)\geq 3$, then the natural map $\H^2_\et(\P^n_k,\mu)\to\H^2_\et(X,\mu)$
is an isomorphism.
\ele
{\em Proof.}
A combination of the weak Lefschetz theorem 
with Poincar\'e duality gives that the map
$\H^i_\et(\P^n_{k_\s},\mu)\to \H^i_\et(X^\s,\mu)$
is an isomorphism for $i<\dim(X)$, see \cite[Cor.~B.6]{K}.
In particular, if $\dim(X)\geq 2$, then
$\H^1_\et(X^\s,\mu)=0$. Then the spectral sequence 
$$E^{p,q}_2=\H^p(k,\H^q_\et(X^\s,\mu))\Rightarrow \H^{p+q}_\et(X,\mu)$$
implies the first claim.

If $\dim(X)\geq 3$, then $\H^2_\et(\P^n_{k_\s},\mu) \to \H^2_\et(X^\s,\mu)$ is an isomorphism
of $\Ga$-modules.
The above spectral sequence gives rise to 
the following commutative diagram with exact rows
$$ \xymatrix{0 \ar[r]& \H^2(k,\mu) \ar[r]& \H^2_\et(X,\mu) \ar[r]& \H^2_\et(X^\s,\mu)^\Ga \ar[r]& 
\H^3(k,\mu)\\
0 \ar[r]& \H^2(k,\mu) \ar[r]\ar[u]^\id& \H^2_\et(\P^n_{k},\mu) \ar[r]\ar[u]& 
\H^2_\et(\P^n_{k_\s},\mu)^\Ga \ar[r]\ar[u]^\cong& \H^3(k,\mu)\ar[u]^\id}    $$
By the 5-lemma we deduce that 
$ \H^2_\et(\P^n_{k},\mu) \to  \H^2_\et(X,\mu)$
is an isomorphism. \hfill$\Box$

\bpr\label{LefschezP4k}
Let $X \subset \P^n_{k}$ be a smooth complete intersection of dimension $\dim(X)\geq 3$.
Then the natural map
$$\H^2(k,\mu) \to \H^2_{\nr}(X,\mu)$$
is an isomorphism.
\epr
{\em Proof.} The map $\Z\cong\Pic(\P^n_{k_{\s}}) \to \Pic(X^\s)$ is an isomorphism by
\cite[Ch.~IV, Cor.~3.2]{H}, hence $\Pic(X) \to \Pic (X^\s)$ is an isomorphism.
 By Proposition \ref{raetsel2} it is thus enough to prove that the map 
$\H^2(k,\mu) \to \H^2_{\nr}(X,\mu)$ is surjective. 

Choose an affine subspace $\A^n_k\subset\P^n_k$ such that $X\cap\A^n_k\not=\emptyset$.
Our map is the composition of maps in the top row of the following natural commutative diagram:
$$\xymatrix{
\H^2(k,\mu) \ar[r]\ar[d]^\id&\H^2_\et(\P^n_k,\mu)\ar[r]^\cong\ar[d]&\H^2_\et(X,\mu) \ar[r]\ar[d]&  
\H^2_\nr(X,\mu)\ar @{^{(}->}[d]\\
\H^2(k,\mu) \ar[r]^{\cong \ \ }&\H^2_\et(\A^n_k,\mu)\ar[r]&\H^2_\et(X\cap\A^n_k,\mu) \ar[r]&  
\H^2(k(X),\mu)}$$
In the top row, the middle map is an isomorphism by Lemma \ref{Lefschetz} (b), 
and the right-hand map is surjective, as was recalled in the introduction. Thus any
$a\in \H^2_\nr(X,\mu)$ can be lifted to an element $b\in \H^2_\et(\P^n_k,\mu)$.
The image of $b$ in $\H^2_\et(\A^n_k,\mu)$ comes from a unique element $c\in \H^2(k,\mu)$.
The commutativity of the diagram gives that the image of $c$ in $\H^2(k(X),\mu)$
is equal to the image of $a$. But the right-hand vertical map is injective, hence $c$
is a desired lifting of $a$ to $\H^2(k,\mu)$. \hfill $\Box$

\section{Generic diagonal hypersurfaces} \label{gen diag}

Let $\Pi_1$ (respectively, $\Pi_2$) be the projective space with homogeneous coordinates
$x_0,\ldots,x_n$ (respectively, $t_0,\ldots,t_n$). Write $K=k(\Pi_2)$.
Let $X\subset \Pi_1\times \Pi_2$ be the hypersurface
\begin{equation}
t_0x_0^d+\ldots+t_nx_n^d=0.\label{e3}
\end{equation}
Let $p$ be the projection $X\to\Pi_1$, and let $f$ be the projection $X\to\Pi_2$.
The generic fibre $X_K$ of $f$
is a smooth diagonal  hypersurface of degree $d$ in the projective space
$(\Pi_1)_K\cong\P^n_K$.

\ble\label{integralfibres}  With notation as above, the following statements hold.

{\rm (i)}  The fibres of $f$ at codimension $1$
points of  $\Pi_{2}$ are integral if $n\geq 2$ and geometrically integral
if $n\geq 3$.

{\rm (ii)} The fibres of $f$ at codimension $2$ points of $\Pi_{2}$
are integral if $n\geq 3$ and geometrically integral if $n\geq 4$.
\ele
{\em Proof.} One only needs to check this for the singular fibres,
which are the fibres above the generic points of the projective subspaces
given by $t_{i}=0$ or by $t_{i}=t_{j}=0$. \hfill $\Box$

\subsection{Unramified cohomology in degree 1}\label{cohodegree1}

\ble \label{rock}
Let $f\colon X\to Y$ be a proper, dominant morphism of smooth and geometrically
integral varieties over a field $k$. 
Write $K=k(Y)$ and let $X_K$ be the generic fibre of $f$. 
Assume that the fibres of $f$ over the points of $Y$ of codimension~$1$ are integral
and $X_K$ is geometrically integral. Let $m \geq 2$ be an integer. Then
the map $f^*\colon \Pic(Y)/m\to\Pic(X)/m$ is injective if and only if 
$\Pic(X)[m]\to \Pic(X_K)[m]$ is surjective.
\ele
{\em Proof.} In our situation we have an exact sequence
\begin{equation}
0\to\Pic(Y)\stackrel{f^*}\lra\Pic(X)\to\Pic(X_K)\to 0.
\label{e2}
\end{equation}
Exactness at $\Pic(X_K)$: since $X$ is smooth, the Zariski closure in $X$
of a Cartier divisor in $X_K$ is a Cartier divisor in $X$. Exactness at $\Pic(X)$:
if $D\in\Div(X)$ restricts to a principal divisor in $X_K$, then $D$ is the sum of a 
principal divisor in $X$ and a divisor contained in the fibres of $f$, which
by our assumption is contained in $f^*\Div(Y)$. Exactness at $\Pic(Y)$: if $D\in\Div(Y)$
is such that $f^*D=\div_X(\phi)$, where $\phi\in k(X)^\times$, then the restriction of $\phi$
to $X_K$ is a regular function. Since $X_K$ proper and geometrically integral, we must have
$\phi\in K^\times$. Then $D-\div_Y(\phi)\in \Div(Y)$ goes to zero in $\Div(X)$, so $D=\div_Y(\phi)$
is a principal divisor in $Y$.

Applying the snake lemma to the commutative diagram obtained from (\ref{e2})
and multiplication by $m$, proves the lemma. \hfill $\Box$

\bpr \label{roll}
Let $m \geq 2$ be an integer. 
Let $k$ be a field of characteristic exponent coprime to $m$.
Let $f\colon X\to Y$ be a proper, dominant morphism of smooth and geometrically
integral varieties over $k$ such that

\smallskip

{\rm (i)} the fibres of $f$ over the points of $Y$ of codimension $1$ are integral
and the generic fibre $X_K$ is geometrically integral (where $K=k(Y)$);

{\rm (ii)} $\Pic(X)$ is torsion-free;

{\rm (iii)} $f^*\colon \Pic(Y)/m\to\Pic(X)/m$ is injective.

\smallskip

\noindent Then $\H^1(K,\mu_m)\to \H^1_\et(X_K,\mu_m)$ is an isomorphism.
\epr
{\em Proof.} The Kummer sequence gives rise to an exact sequence
$$0\to K^\times/K^{\times m}\to \H^1_\et(X_K,\mu_m)\to \Pic(X_K)[m]\to 0.$$
By Lemma \ref{rock} we have $\Pic(X_K)[m]=0$. \hfill $\Box$

\bthe \label{deg1bis}
Let $n\geq 2$. Let $\Pi_{1}$, $\Pi_{2}$, $X$, $K=k(\Pi_2)$ be as above.
Then the map $\H^1(K,\mu) \to \H^1_\et(X_{K},\mu)$ is an isomorphism.
\ethe
{\em Proof.}
Let us first prove the statement for $\mu=\mu_{m}$ with
$m$ not divisible by ${\rm char}(k)$. 
Let us check the assumptions of Proposition \ref{roll}
for $f \colon X \to \Pi_{2}$.
By Lemma \ref{integralfibres}, assumption  (i)
is satisfied. The projection $p\colon X \to \Pi_{1}$ is a projective  bundle over  
$\Pi_{1}$. We have a commutative diagram with exact rows
$$\xymatrix{
 0 \ar[r]&\Pic(\Pi_{1}) \ar[r]& \Pic(X) \ar[r]& \Pic(\P^{n-1}_{k(\Pi_1)}) \ar[r]& 0\\
 0 \ar[r]&\Pic(\Pi_{1}) \ar[r]\ar[u]^\id& \Pic(\Pi_1\times\Pi_2) \ar[r]\ar[u]& 
\Pic((\Pi_2)_{k(\Pi_1)}) \ar[r]\ar[u]^\cong& 0}$$
The right-hand vertical map is induced by 
the inclusion of a projective hyperplane in a projective space, so it is an isomorphism.
Hence (ii) holds and the restriction map
$\Pic(\Pi_1\times \Pi_2)\to\Pic(X)$ is an isomorphism. It follows that 
$\Pic(\Pi_2)\to\Pic(X)$ is split injective, hence (iii) holds.

For an arbitrary group $\mu$,
let $E/k$ be a finite Galois extension, with Galois group $G$,
such that $\mu_{E}= \mu\times_{k}E$ is isomorphic to a finite product of groups $\mu_{m,E}$
where $m$ is coprime to ${\rm char}(k)$.
Let $L$ be the compositum of the linearly disjoint field extensions $K/k$ and $E/k$.
We have $\mu(E)=\mu(L)= \H^0_\et(X_{L},\mu)$. The Hochschild--Serre spectral sequence gives rise to
the following commutative diagram with exact rows
$$\xymatrix{0 \ar[r]& \H^1(G, \mu(L)) \ar[r]&  \H^1_\et(X_{K},\mu) \ar[r]&  \H^1_\et(X_{L}, \mu)^G \ar[r]& \H^2(G, \mu(L))\\
0 \ar[r]& \H^1(G,\mu(L)) \ar[r]\ar[u]^\id& \H^1(K, \mu) \ar[r]\ar[u]& 
\H^1(L,\mu)^G \ar[r]\ar[u]^\cong& \H^2(G,\mu(L))\ar[u]^\id}$$
Since the result is already proved for $\mu_m$, all vertical maps, except possibly the map
 $\H^1(K, \mu)  \to \H^1_\et(X_{K}, \mu)$, are isomorphisms. Hence so is
 this map.
\hfill $\Box$

\begin{rem} 
{\rm
The geometric argument 
based on the projective bundle structure of
 $X \subset \Pi_{1} \times \Pi_{2}$ over $\Pi_{1}$  in 
the proof of Theorem \ref{deg1bis}
is needed  only in the case $n=2$, that is, when the hypersurface $X_{K} \subset \P^2_{K}$ is a smooth curve of degree $d$.
When $n\geq 3$ and $X \subset \P^n_{K}$ is an {\it arbitrary} smooth
hypersurface, we have
$\H^1(K,\mu) \cong \H^1(X_{K},\mu)$ by Lemma \ref{Lefschetz} (a).}
\end{rem}

 \subsection{Basic diagram}\label{basicdiagram}

{\em We now assume $n\geq 3$ and $i \geq 2$.}
Recall the Bloch--Ogus complex from Section \ref{unramifcoho}:
$$\H^i(k(X),\mu )\stackrel{(\partial_x)\,}\lra \bigoplus_{x\in X^{(1)}}\H^{i-1}(k(x),\mu(-1))\to 
\bigoplus_{x\in X^{(2)}}\H^{i-2}(k(x),\mu(-2)).$$
Since the fibres $X_y=f^{-1}(y)$ over $y\in \Pi_2^{(1)}$ 
are integral (which holds for $n\geq 2$, see Lemma \ref{integralfibres}) 
we obtain a complex
$$\H^i_\nr(X_K,\mu)\stackrel{(\partial_y)}\lra 
\bigoplus_{y\in \Pi_2^{(1)}}\H^{i-1}(k(X_y),\mu(-1))\to 
\bigoplus_{x\in X^{(2)}}\H^{i-2}(k(x),\mu(-2)).$$
To simplify notation, in what follows we do not write the coefficients of cohomology groups.
One should bear in mind that there is a change of twist when the codimension of points increases.

Since this is a complex,
the image of $\partial_y$ is unramified over the smooth locus of $X_y$.
If $X_y$ is smooth we write $X_y'=X_y$. In the opposite case,
$X_y$ is the projective cone over
the hyperplane section of $X$ given by some 
$t_{i}=0$,
and then we denote by $X_y'$
this hyperplane section, 
which is smooth since $n\geq 3$.
In this case, the smooth locus $X_{y, {\rm sm}}\subset X_y$ is an affine bundle over $X_y'$, so
we have $\H^{i-1}_\nr(X_{y, {\rm sm}})\cong \H^{i-1}_\nr(X_y')$ by (\ref{affine_bundle}).
Thus $\Im(\partial_y)$ is contained in $\H^{i-1}_\nr(X_y')$.
Since the fibres $X_y$ over $y\in \Pi_2^{(2)}$ 
are integral (note that they need not be
 geometrically integral if $n=3$), from the diagram in Section \ref{unramifcoho} 
we obtain a commutative diagram of complexes
$$\xymatrix{0\ar[r]&\H^i_\nr(X_K)/\H^i(k)\ar[r]&\bigoplus_{y\in \Pi_2^{(1)}}
\H^{i-1}_\nr(X_y')\ar[r]&
\bigoplus_{y\in \Pi_2^{(2)}}\H^{i-2}(k(X_y))\\0\ar[r]&
\H^i(K)/\H^i(k)\ar[r]\ar[u]&\bigoplus_{y\in \Pi_2^{(1)}}\H^{i-1}(k(y))\ar[r]\ar[u]&
\bigoplus_{y\in \Pi_2^{(2)}}\H^{i-2}(k(y))\ar[u]}$$
where the vertical maps are induced by $f$. Note that 
since $X$ is a projective bundle over the projective space $\Pi_{1}$,
the maps $\H^{i}(k) \to \H^{i}(k(X))$ is injective.  So is the map
 $\H^{i}(k) \to \H^{i}(K)=  \H^{i}(k(\Pi_{2}))$.

Let $Y=\A^n_k\subset\Pi_2$ be the affine space given by $t_0\neq 0$. 
From the previous diagram we then get a commutative diagram of complexes
\begin{equation}
\begin{split}
\xymatrix{0\ar[r]&\H^i_\nr(X_K)/\H^i(k)\ar[r]&\bigoplus_{y\in Y^{(1)}}\H^{i-1}_\nr(X_y')\ar[r]&
\bigoplus_{y\in Y^{(2)}}\H^{i-2}(k(X_y))\\0\ar[r]&
\H^i(K)/\H^i(k)\ar[r]\ar[u]&\bigoplus_{y\in Y^{(1)}}\H^{i-1}(k(y))\ar[r]\ar[u]&
\bigoplus_{y\in Y^{(2)}}\H^{i-2}(k(y))\ar[u]} \label{diag}
\end{split}
\end{equation}
Since $Y\cong\A^n_k$, the bottom complex is exact by \cite[Prop.~8.6]{R96}.

The homology group of the top complex at the first term is $\H^i_\nr(X_Y)/\H^i(k)$, where
$X_Y=f^{-1}(Y) \subset X$. 
Let us show that this group is zero. The fibres of  $p\colon X\to \Pi_1$ are hyperplanes in $\Pi_2$.
The map  $p\colon X_Y\to U$ is an affine bundle, and
$p(X_Y)=U$, where $U=\P^n_k\setminus \{(1:0:\ldots:0)\}$.
By (\ref{affine_bundle}) the map $p^*\colon \H^i_\nr(U)\to\H^i_\nr(X_Y)$
is an isomorphism. 
Since $U$ is the complement to a $k$-point in 
$\Pi_1\cong\P^n_k$, and $n\geq 2$, we have
$$\H^{i}(k,\mu) \cong \H^{i}_{\nr}(\Pi_{1},\mu) \cong \H^{i}_{\nr}(U, \mu).$$

The following lemma is proved by a straightforward diagram chase.

\ble \label{basic}
Suppose that we have a commutative diagram of abelian groups
$$\xymatrix{&A \ar@{^{(}->}[r]^i& B\ar[r]^j& C\\
0\ar[r]&D\ar[r]\ar[u]^a&E\ar[r]\ar[u]^b_\cong
&F\ar@{_{(}->}[u]^c
}
$$
where $i$ is injective, $b$ is an isomorphism, $c$ is injective, the top row is a complex,
and the bottom row is exact. Then $a$ is an isomorphism.
\ele

From Lemma \ref{basic} we conclude:

\begin{prop}\label{mecanisme}
With notation as above, if the middle vertical map in diagram $(\ref{diag})$
is an isomorphism and the right-hand vertical map is injective, then
$f^*\colon\H^i(K,\mu)\to\H^i_\nr(X_K,\mu)$ is an isomorphism.
\end{prop}

\subsection{Unramified cohomology in degrees 2 and 3}\label{cohodegree2and3}

  The main result of this paper is the following

\bthe \label{main}
Let $\Pi_1$ (respectively, $\Pi_2$) be the projective space with homogeneous coordinates
$x_0,\ldots,x_n$ (respectively, $t_0,\ldots,t_n$). Write $K=k(\Pi_2)$.
Let $X\subset \Pi_1\times \Pi_2$ be the hypersurface
\begin{equation}
t_0x_0^d+\ldots+t_nx_n^d=0.\label{e4}
\end{equation}
Let $f\colon X\to\Pi_2$ be the natural projection, and let $X_K$ be the generic fibre of $f$.
Let $\mu$ be a finite \'etale commutative group $k$-scheme of order not divisible by ${\rm char}(k)$.

\smallskip

{\rm (i)} If $n\geq 3$, then $f^*\colon\H^2(K, \mu)\to 
\H^2_\nr(X_K,\mu)$ is an isomorphism.

{\rm (ii)} If $n\geq 4$, then $f^*\colon\H^3(K,\mu) \to \H^3_\nr(X_K,\mu)$ is an isomorphism.

\ethe
{\em Proof.} (i) Consider diagram (\ref{diag}) for $i=2$. 
Then the middle vertical map  of the diagram
is an isomorphism. This follows from Theorem \ref{deg1bis} 
when
$X_y$ is singular, which happens exactly when the codimension 1 point  $y$ is given by $t_i=0$ for some $i=1,\ldots, n$.
(Note that if $n=3$  we then need   Theorem \ref{deg1bis} in the case $n=2$.)
If $X_y$ is smooth, the isomorphism
follows from Lemma \ref{Lefschetz} (a).
 By Lemma \ref{integralfibres}, each fibre $X_{y}$ at a codimension 2 point $y$  is integral,
 hence the  right hand vertical map is
 injective.
 By Proposition \ref{mecanisme}, this proves (i).

(ii) Consider  diagram (\ref{diag}) for $i=3$. 
For $y\in Y^{(1)}$ such that $X_y$ is singular, the vertical map 
$\H^{2}(k(y)),\mu(-1)) \to \H^{2}_\nr(X_y',\mu(-1))$
is an isomorphism by (i).  For 
$y\in Y^{(1)}$ such that $X_{y}$ is smooth, 
the map $ \H^{2}(k(y),\mu(-1)) \to \H^{2}_\nr(X_y,\mu(-1))$
is an isomorphism by Proposition  \ref{LefschezP4k}.
For $y\in \Pi_2^{(2)}$
the fibre $X_{y}$ is geometrically integral over $k(y)$ by Lemma \ref{integralfibres},
hence $k(y)$ is separably closed in $k(X_y)$. Thus the restriction map
$\H^1(k(y), \mu(-2)) \to  \H^1(k(X_{y}),\mu(-2))$
is injective, so the right-hand vertical map in the diagram is  injective. 
By Proposition \ref{mecanisme}, this proves (ii). \hfill $\Box$

\bco \label{Br}
For $n \geq3$, the map $\Br(K)\to\Br(X_K)$ induces an isomorphism of subgroups of elements
of order not divisible by ${\rm char}(k)$.
\eco
{\em Proof.} 
This follows from  Theorem \ref{main} (i) by taking $\mu=\mu_{m}$
for each  integer $m$ not divisible by ${\rm char}(k)$.
\hfill $\Box$

 \begin{rem}{\rm
 Only the case $n=3$  of this corollary requires the above proof. For $n\geq 4$ and any
 smooth hypersurface in $\P^n$, we have the general  Proposition \ref{LefschezP4k}.
 }
\end{rem}

\section{Pairs of polynomials} \label{pairs}

In this section we give a short elementary proof that the Brauer group of the surface
given by the equation (\ref{pair}) below 
over the field of rational functions $K=k(t)$, with $t=\lambda/\mu$,
is naturally isomorphic to $\Br(K)$ (away from $p$-primary torsion if ${\rm char}(k)=p$).
In the case when $k$ has characteristic zero, this follows from more general results of 
\cite{GS}, namely, the combination of \cite[Thm.~1.1 (i)]{GS} and \cite[Thm.~1.4]{GS}.

\bthe \label{f=ag}
Let $k$ be a field. Let $d$ be a positive integer.
Let $f(x,y)$ and $g(z,t)$ be products of $d$ pairwise non-proportional linear forms.
Let $X \subset  \P^1_{k} \times_{k}\P^3_{k}$ be the hypersurface given by 
\begin{equation}
\lambda f(x,y)= \mu g(z,t), \label{pair}
\end{equation}
where $(\lambda:\mu)$ are homogeneous coordinates in $\P^1_{k}$ and $(x:y:z:t)$ are
homogeneous coordinates in $\P^3_{k}$.
Let $K=k(\P^1_k)$ and let $X_{K}$ be the generic fibre of the projection $f\colon X \to \P^1_k$.
Then the natural map $\Br(K)\to\Br(X_{K})$ induces an isomorphism of subgroups of
elements of order not divisible by ${\rm char}(k)$.
\ethe
{\em Proof.} 
The singular locus $X_{\rm sing}$ is contained in the union of fibres of $f$ above $\lambda=0$ and $\mu=0$.
The fibre above $\mu=0$ is given by $f(x,y)=0$. It is a union of
$d$ planes in  $\P^3_k$ through the line $x=y=0$.
The intersection of $X_{\rm sing}$ with the fibre above $\mu=0$ is the zero-dimensional scheme  
 given by $x=y=g(z,t)=0$.  The situation above $\lambda=0$ is entirely similar.
 Let $Y=X\setminus X_{\rm sing}$ be the smooth locus of $X/k$.
The projection $p\colon X \to \P^3_{k}$ is a birational morphism 
which restricts to an isomorphism $Y_{V} \tilde\lra V$ 
on the complement $V$ to the curve in $\P^3_k$ given by $f(x,y)=g(z,t)=0$. 
We have
$$\Br(k) \cong \Br(\P^3_{k}) \cong \Br(V) \cong \Br(Y_{V}),$$ 
where the first isomorphism is by \cite[Thm.~6.1.3]{CTS21} and the second one
is by purity for the Brauer group \cite[Thm.~3.7.6]{CTS21}. Since $Y(k)\not=\emptyset$, we have
$\Br(k)\subset\Br(Y)\subset\Br(Y_V)$ where the second inclusion is
by \cite[Thm.~3.5.5]{CTS21}. We conclude that $\Br(Y)\cong\Br(k)$.

Let $m\geq 2$ be an integer not divisible by ${\rm char}(k)$. If a closed fibre
$X_{M}=f^{-1}(M)$ is smooth, then $X_M$ is a
smooth surface in $\P^3_{k(M)}$, thus we have
\begin{equation}
\H^1_\et(X_{M}, \Z/m)\cong\H^1(k(M),\Z/m)\label{raz}
\end{equation} 
by Lemma \ref{Lefschetz} (a).
The smooth locus of the fibre of $f$  above $\mu=0$ 
is a disjoint union of $d$ affine planes $\A^2_k$.
We have 
\begin{equation}
\H^1_\et(\A^2_{k},\Z/m)\cong\H^1(k,\Z/m) \label{dva}
\end{equation} since ${\rm char}(k)$ does not divide $m$.

Without loss of generality we can write 
$$f(x,y)=c\prod_{i=1}^d(x-\xi_i y), \quad\quad
g(z,t)=c'\prod_{j=1}^d(z-\rho_j t),$$
where $c,c'\in k^\times$ and $\xi_i,\rho_j\in k$ for $i,j=1,\dots,d$. 
We note that for each pair $(i,j)$ the map 
$s_{ij}\colon(\lambda:\mu)\to\big((\lambda:\mu), (\xi_i:1:\rho_j:1)\big)$ 
is a section of the morphism $f\colon X \to \P^1_{k}$.

Each section $s_{ij}$ gives a $K$-point of $X_K$. Thus the natural map $\Br(K)\to\Br(X_K)$ is injective.

Let $\alpha\in \Br(X_{K})[m]$. Evaluating $\alpha$ at the $K$-point of $X_K$ given by $s_{1,1}$
gives an element $\beta\in\Br(K)[m]$. We replace $\alpha$ by $\alpha-\beta$. 

Note that each section $s_{ij}(\P^1_k)$ meets every closed fibre of $f$ at a smooth point.
The new element $\alpha\in\Br(X_K)[m]$ has trivial residue on the
irreducible component of the smooth locus of every fibre of $f$ that $s_{1,1}(\P^1_k)$  intersects. Indeed, by (\ref{raz}) and (\ref{dva})
this residue is constant, but specialises to zero at the intersection point with $s_{1,1}(\P^1_k)$.
In particular, $\alpha$ has trivial residues at the smooth fibres of 
$f$, as well as at the affine plane given by $x-\xi_1 y=0$ in the fibre $\mu=0$
and the affine plane given by $z-\rho_1 t=0$ in the fibre $\lambda=0$.

We now evaluate $\alpha$ at the $K$-point of $X_K$ given by $s_{1,j}$, where $j=2,\ldots,d$.
The result is an element of $\Br(K)$ which is unramified everywhere except possibly
at the $k$-point of $\P^1_k$ given by $\lambda=0$. 
By Faddeev reciprocity, 
the residue at that point must be zero, too.
This implies that $\alpha$ is unramified at the smooth locus of the fibre at $\lambda=0$.
A similar argument using sections $s_{i,1}$ for $i=2,\ldots,d$ shows that
$\alpha$ is unramified at the smooth locus of the fibre at $\mu=0$.

We see that the residue of $\alpha$ at every point of  codimension 1
of $Y$ is zero. Thus $\alpha$ belongs to 
$\Br(Y)$, hence to $\Br(k)$. 
We conclude that $\Br(K)[m]\stackrel{\sim}\lra\Br(X_K)[m]$. \hfill $\Box$

\noindent Universit\'e Paris-Saclay, CNRS, Laboratoire de math\'ematiques d’Orsay, 91405, Orsay, 
France.

\medskip

\noindent \texttt{jean-louis.colliot-thelene@universite-paris-saclay.fr}

\bigskip

\noindent Department of Mathematics, South Kensington Campus,
Imperial College London, SW7 2BZ England, U.K. -- and --
Institute for the Information Transmission Problems,
Russian Academy of Sciences, 19 Bolshoi Karetnyi, Moscow, 127994
Russia

\medskip

\noindent \texttt{a.skorobogatov@imperial.ac.uk}

\end{document}